\documentclass[12pt]{article}

\usepackage{amsfonts,amsmath,amsthm}

\usepackage[english,francais]{babel}

\newtheorem{thm}{Theorem}
\newtheorem{lem}[thm]{Lemma}
\newtheorem{prop}[thm]{Proposition}
\newtheorem{cor}[thm]{Corollary}
\newtheorem{remark}{Remark}
\newtheorem{example}{Example}

\def \ZZ {\mathbb{Z}}
\def \QQ {\mathbb{Q}}

\def \EE {\mathbb{E}}
\def \NN {\mathbb{N}}
\def \RR {\mathbb{R}}
\def \PP {\mathbb{P}}
\def \kL{\mathcal{L}}
\def \kA{\mathcal{A}}
\title{On Wiener-Hopf factors for stable processes}
\author{Piotr Graczyk \footnote{This research was partially supported by grants MNiSW N N201 373136  and ANR-09-BLAN-0084-01.}\\
\scriptsize  LAREMA, Universit\'e d'Angers, 
2 Bd Lavoisier, \\
 \scriptsize 49045 Angers Cedex 1, France\\
 \scriptsize Piotr.Graczyk@univ-angers.fr\\\\
Tomasz Jakubowski \footnote{This research was partially supported by grants MNiSW 3971/B/H03/2009/37, ANR-09-BLAN-0084-01 and the fellowship of CCRRDT Pays de la Loire.}\\
\scriptsize LAREMA, Universit\'e d'Angers, 
2 Bd Lavoisier, \\
 \scriptsize 49045 Angers Cedex 1, France\\
\scriptsize Institute of Mathematics, Wroc\l{}aw University of Technology, Wyb. Wyspia\'nskiego 27,\\
\scriptsize 50-370 Wroc\l{}aw, Poland\\
\scriptsize Tomasz.Jakubowski@pwr.wroc.pl \tiny}
\date{\empty}
\begin{document}
\maketitle
\selectlanguage{english}
\begin{abstract}
We give a series representation of the logarithm of the bivariate Laplace exponent $\kappa$ of $\alpha$-stable processes for almost all $\alpha \in (0,2]$.
\end{abstract}
\selectlanguage{francais}
\begin{abstract}
Nous donnons un d\'eveloppement en s\'erie
du logarithme de l'exposant de Laplace bivari\'e
$\kappa$ des processus $\alpha$-stables
pour presque tous $\alpha\in (0,2]$.
\end{abstract}
\selectlanguage{english} 



\section{Introduction}
The fluctuation theory of L\'evy processes is one of
 the domains of probability very
 actively
developing in the  last years, and with important applications
in mathematical finance; cf. the recent monograph of A. Kyprianou \cite{MR2250061}
and papers \cite{MR2440923}, \cite{MR2451576}, \cite{MR2499867}, \cite{MR2483728}. The $\alpha$-stable L\'evy processes play a 
primordial role in this theory. We address in this article
one of the key  problems of the Wiener-Hopf factorization theory
of the $\alpha$-stable  processes: the computation
 of the bivariate Laplace exponent $\kappa(\gamma,\beta)$.

The aim of this paper is to give a series representation of the integral 
\begin{equation}\label{eq:int}
g(\beta) = \frac{\sin(\pi \rho)}{\pi}\int_0^\infty \frac{\beta \log(1 + x^\alpha)}{x^2+2x\beta\cos(\pi\rho)+\beta^2}\,dx
\end{equation}
for almost all $\alpha \in (0,2]$ and $\rho \in [1-1/\alpha,1/\alpha] \cap (0,1)$. This integral plays an important role in the theory of stable processes. Using (\ref{eq:int}) one may express the bivariate Laplace exponent $\kappa(\gamma,\beta)$ of the ascending ladder process built from the $\alpha$-stable process $X_t$ with index of stability $\alpha$ and $\rho= \PP(X_1 >0)$ (see e.g. \cite{MR1406564}, \cite{MR2250061}). Namely
\begin{align*}
\gamma^\rho \exp\left\{g(\beta \gamma^{-1/\alpha})\right\} &= \kappa(\gamma,\beta)\\
&= k \exp\left\{\int_0^\infty\int_{(0,\infty)}\frac{e^{-t}-e^{-\gamma t}e^{-\beta x}}{t} \PP(X_t \in dx)\,dt\right\}\,.   
\end{align*}

The integral (\ref{eq:int}) was introduced by Darling in \cite{MR0080393} for $\rho=1/2$ and calculated in the case $\alpha=1$ and $\rho=1/2$, which corresponds to the symmetric Cauchy process and later by Bingham \cite{MR0415780} for spectrally negative stable processes ($1/\rho = \alpha \in (1,2)$). Doney in \cite{MR905336} calculated it for the set of parameters  $(\alpha,\rho)$ satisfying $\rho +k = l/\alpha$ for some $k \in \NN \cup \{0\}$ and $l \in \NN$.  Although the function $\kappa$ plays an important role in the theory of stable (in general L\'evy) processes, the only known closed expression for it is due to Doney. In this note we expand the function $g$ to a power series for almost all $\alpha$ and $\rho$. We denote by $\kL$ the set of Liouville numbers, which will be defined in Section 2. Let $\kA = (0,2] \setminus (\QQ \cup \kL)$. We note that if $\alpha \in \kA$ then by Lemma \ref{Liouv}, $1/\alpha \in \kA$. The main result of this paper is
\begin{thm}\label{thm:main}
Let $\alpha \in \mathcal{A}$, $\rho \in [1-1/\alpha,1/\alpha] \cap (0,1)$ and $0<\beta<1$.  Then
\begin{align}\label{eq:thm}
g(\beta)
& =  \sum_{m=1}^\infty \frac{(-1)^{m+1} \beta^{m}\sin(\rho m \pi)}{m\sin(\frac{m\pi}{\alpha})} + \sum_{k=1}^\infty \frac{(-1)^{k+1} \beta^{\alpha k}\sin(\rho\alpha k\pi)}{k\sin(\alpha k\pi)}\,.
\end{align}
\end{thm}
We note that in view of Lemma \ref{lem:1/beta} it suffices to consider only $0<\beta<1$.  
The series unfortunately does not converge for irrational numbers $\alpha \in \mathcal{L} \cap (0,2)$, but $\mathcal{L}$ has a Lebesgue measure $0$ hence $\mathcal{A}$ contains almost all $\alpha \in (0,2)$. We obtain also a formula for rational $\alpha$ (Proposition \ref{lem:rational}) but the expression is not so closed as in Theorem \ref{thm:main}. If $\rho + k = l/\alpha$ for some integers $k$ and $l$ the formula (\ref{eq:thm}) may be simplified, in particular one may obtain results achieved by Doney in \cite{MR905336} (see Remark \ref{rem:1}). 

The formula (\ref{eq:thm}) opens a way to applications for the study of various
functionals of an $\alpha$-stable  L\'evy process, in particular of the long time behavior
of the supremum process or the law of the  first passage time, cf. the recent results
of Bernyk, Dalang and Peskir \cite{MR2440923} and Kuznetsov \cite{Kuz}. We also profit
from the Theorem 1 in a forthcoming work \cite{PGTJ}, devoted to the first passage time
of symmetric stable processes. 

The paper is organized as follows. In Section 2 we define Liouville numbers and prove some auxiliary lemmas. In Section 3 we prove the main Theorem \ref{thm:main}. In Section 4 we give some remarks, applications and examples.

\section{Liouville numbers}
A number $x \in \RR$ is called a Liouville number if it may be well approximated by rational numbers. More precisely for any $n \in \NN$ there exist infinitely many pairs of integers $p,q$ such that (see e.g. \cite{MR781734})
$$
0 < \left|x - \frac{p}{q}\right| < \frac{1}{q^n}\,.
$$
We denote by $\mathcal{L}$ the set of all Liouville numbers. First we note
\begin{lem}\label{Liouv}
$x \in \mathcal{L}$ if and only if $1/x \in \mathcal{L}$.
\end{lem}
The proof does not seem available in the literature. The following proof was proposed by M. Waldschmidt (\cite{Waldschmidt}).
\begin{proof}
Let $x \not\in \mathcal{L}$. There are $c \in \RR$ and $d \in \NN$ such that for all $p\in \ZZ ,q \in \NN$,
$$
\left|x - \frac{p}{q}\right| \ge \frac{c}{q^d}\,.
$$
Let $p\in \ZZ ,q \in \NN$. We may and do suppose that $|1/x - p/q| <1$. Hence $|p|/q < (|x|+1)/|x|$ and 
$$
\left|\frac{1}{x} - \frac{p}{q}\right| = \frac{|p|}{q|x|}\left|x - \frac{q}{p}\right| \ge \frac{|p|}{q|x|}\, \frac{c}{|p|^d}  \ge \frac{|x|^{d-2}}{(1+|x|)^{d-1}}\, \frac{c}{q^d}\,.
$$
\end{proof}

\begin{lem}\label{lem:Liouv1}
For any $x \in \kA$ and $\beta \in (0,1)$ we have
$$
\sum_{m=1}^\infty \frac{\beta^m}{|\sin(m x\pi)|} < \infty\,.
$$
\end{lem}
\begin{proof}
Since  $x \in \RR \setminus (\QQ \cup \mathcal{L})$, there is $N \in \NN$ such that $\left|x  - \frac{p}{q}\right| > \frac{1}{q^N}$ for all integers $p,q>0$. Hence $|\sin(m x\pi)| > \frac{1}{2m^{N-1}}$ and the lemma follows.
\end{proof}

In the sequel we will need following formulas taken from \cite{{MR2360010}} (formulas 1.445.7, 1.422.3, 1.353.1)
\begin{equation}\label{eq:tmp1}
\sum_{m=1}^\infty (-1)^{m+1} \frac{m\sin(m z)}{m^2 - w^2} = \frac{\pi}{2}\frac{\sin(zw)}{\sin(w\pi)}\,, \qquad z \in (-\pi,\pi), w \in \RR \setminus \ZZ\,,
\end{equation}

\begin{equation}\label{eq:tmp2}
\frac{\pi}{\sin(\pi z)} = \frac{1}{z} - \sum_{k=1}^\infty \frac{(-1)^k 2 z}{k^2 - z^2}, \qquad z \in \RR \setminus \ZZ\,,
\end{equation}

\begin{equation}\label{eq:tmp3}
\sum_{k=1}^{n-1} p^k \sin(kx) = \frac{p\sin(x) - p^n\sin(nx)+p^{n+1}\sin((n-1)x)}{1-2p\cos(x)+p^2}\,.
\end{equation}

\begin{lem}\label{lem:Liouv2}
Let $\alpha \in \kA$ and $\rho \in [1-1/\alpha,1/\alpha] \cap (0,1)$. Then there are constants $C$ and $N$ such that for all $M,k \in \NN$
$$
\sum_{m=1}^M \frac{(-1)^m \sin(m \rho \pi) m}{m^2 - (\alpha k)^2} \le C k^N\,.
$$
\end{lem}
\begin{proof}
Let $K$ be the smallest integer larger then $\alpha k +1$. Like in the proof of Lemma \ref{lem:Liouv1} we take $N$ such that $\left|\alpha   - \frac{p}{q}\right| > \frac{1}{q^N}$ for all integers $p,q$. Then $|m^2- (\alpha k)^2| > mk^{-N+1}$ for all $m,k \in \NN$ and we get
$$
\left|\sum_{m=1}^{K-1} \frac{(-1)^m \sin(m \rho \pi) m}{m^2 - (\alpha k)^2}\right| \le (\alpha k +1)k^{N-1} \le 3k^N\,.
$$
Denote $a_m = (-1)^m \sin(m \rho \pi) $ and $b_m = \frac{m}{m^2 - (\alpha k)^2}$. By (\ref{eq:tmp3}) for any $M \ge 1$ we have
$$
\left|\sum_{m=1}^M a_m \right| \le \frac{3}{2(1+\cos(\rho\pi))} = c\,.
$$
 Since $b_m$ is decreasing for $m \ge K$ we get 
\begin{align*}
\left|\sum_{m=K}^M a_m b_m\right| &=\left| \sum_{m=K}^{M-1}  (b_m - b_{m+1}) \sum_{n=K}^m a_n + b_M  \sum_{n=K}^M a_n \right| \\
& \le  \sum_{m=K}^{M-1}  (b_m - b_{m+1}) \left|\sum_{n=K}^m a_n\right| + b_M  \left|\sum_{n=K}^M a_n \right| \le 2c b_K\ \le 2c.
\end{align*}
\end{proof}

\section{Proof of Theorem \ref{thm:main}}
The following lemma justifies our restriction in Theorem \ref{thm:main} to $0 < \beta <1$
\begin{lem}\label{lem:1/beta}
$$\kappa(1,\beta) = \kappa(1,1/\beta) \beta^{\alpha\rho}\,.$$
\end{lem}
\begin{proof}
After substituting $x=1/y$ we get
\begin{align*}
g(\beta) &= \beta \frac{\sin(\pi \rho)}{\pi}\int_0^\infty \frac{\log(1 + y^\alpha) - \log(y^\alpha)}{1+2y\beta\cos(\pi\rho)+y^2\beta^2}\,dy \\
&= g(1/\beta) + \frac{\alpha \sin(\pi \rho)}{\pi}\int_0^\infty \frac{\log(\beta) - \log(z)}{1+2z\cos(\pi\rho)+z^2}\,dz \\
& = g(1/\beta) + \frac{\alpha \log(\beta)}{\pi}\int_0^\infty \frac{\sin(\pi \rho)}{1+2z\cos(\pi\rho)+z^2}\,dz = g(1/\beta) + \alpha\rho \log(\beta)
\end{align*}
and the lemma follows.
\end{proof}
A derivative of the function $g$  is equal to
\begin{align*}
g'(\beta) &= \frac{\partial}{\partial \beta} \frac{\sin(\pi \rho)}{\pi}\int_0^\infty \frac{\log(1 + \beta^\alpha x^\alpha)}{x^2+2x\cos(\pi\rho)+1}\,dx \\
& =\frac{\sin(\pi\rho)\alpha}{\pi}\int_0^\infty \frac{x^\alpha}{1+x^\alpha}\frac{1}{x^2+2x\beta\cos(\pi\rho)+\beta^2}\,dx\,.
\end{align*}
Our aim is to prove
\begin{lem}\label{lem:main}
Let $\alpha \in \mathcal{A}$, $\rho \in [1-1/\alpha,1/\alpha] \cap (0,1)$ and $0<\beta<1$. 
Then
\begin{align*}
& g'(\beta)  =  \sum_{m=1}^\infty \frac{(-1)^{m+1} \beta^{m-1}\sin(\rho m\pi)}{\sin(\frac{m\pi}{\alpha})} + \alpha\sum_{k=1}^\infty \frac{(-1)^{k+1} \beta^{\alpha k-1}\sin(\rho\alpha k \pi)}{\sin(\alpha k\pi)}\,.
\end{align*}
\end{lem}

\begin{lem}\label{lem:int0b}
For any $p>0$ and $0<b<1$
\begin{equation*}
\int_0^b \frac{y^{p}}{1+y}dy =  \sum_{k =0}^\infty \frac{(-1)^{k}b^{k+1+p}}{k+1+p}\,.
\end{equation*}
\end{lem}
\begin{proof}
By Fubini theorem
\begin{eqnarray*}
\int_0^b \frac{y^{p}}{1+y}dy = \int_0^b\sum_{k=0}^\infty (-1)^k y^{p+k} dy = 
 \sum_{k =0}^\infty \frac{(-1)^{k}b^{k+1+p}}{k+1+p}\,.
\end{eqnarray*}
\end{proof}

\begin{lem}\label{lem:intbinfty}
For any $0<b\le 1$ and $p \in (0,\infty) \setminus \NN$ we have
\begin{equation}\label{eq:intbtoinfty}
\int_b^\infty \frac{y^{-p}}{1+y}dy = \frac{\pi}{\sin(p\pi)} + \sum_{k=0}^{\infty} \frac{(-1)^{k+1}b^{k+1-p}}{k+1-p}\,.
\end{equation}
\end{lem}
\begin{proof}
Since the derivatives in $b$ of both sides of (\ref{eq:intbtoinfty}) are equal we have for $b \in (0,1)$
$$
\int_b^\infty \frac{y^{-p}}{1+y}dy = C + \sum_{k=0}^{\infty} \frac{(-1)^{k+1}b^{k+1-p}}{k+1-p}\,.
$$
To calculate the constant $C$ we take $b \to 1$ and by (\ref{eq:tmp2}) we get 
\begin{align*}
C &= \int_1^\infty \frac{y^{-p}}{1+y}dy - \sum_{k=0}^{\infty} \frac{(-1)^{k+1}}{k+1-p}  = \int_0^1 \frac{x^{p-1}}{1+x}dx - \sum_{k=0}^{\infty} \frac{(-1)^{k+1}}{k+1-p}\\
& = \sum_{k=0}^\infty (-1)^k \int_0^1 x^{p+k-1}dx - \sum_{k=0}^{\infty} \frac{(-1)^{k+1}}{k+1-p} = \sum_{k=0}^\infty  \frac{(-1)^k}{k+p} + \sum_{k=0}^{\infty} \frac{(-1)^k }{k+1-p}\\
& = \frac{1}{p} - \sum_{k=1}^{\infty} (-1)^k \frac{ 2p}{k^2-p^2} = \frac{\pi}{\sin(p\pi)}\,.
\end{align*}
\end{proof}
Since for any $n \in \NN$ 
$$
\lim_{p \to n} \left(\frac{\pi}{\sin(p\pi)} + \frac{(-1)^{n}b^{n-p}}{n-p}\right) = (-1)^n \ln b\,,
$$
we get 
\begin{cor}\label{cor:rational}
For $p \in \NN$ and $0<b<1$
$$
\int_b^\infty \frac{y^{-p}}{1+y}dy = (-1)^p \ln b + \sum_{k \in \NN, k\not=p-1} \frac{(-1)^{k+1}b^{k+1-p}}{k+1-p}\,.
$$
\end{cor}
\begin{proof}[Proof of Lemma \ref{lem:main}]
We note that (see \cite[1.447.1]{MR2360010})
\begin{align}
\sum_{m=0}^\infty (-1)^{m} x^m\sin((m+1)z) = \frac{\sin(z)}{x^2+2x\cos(z) +1}\,, \qquad |x|<1\,.  \label{eq: sinseries1}
\end{align}
First we will calculate $\int_0^\beta$. From (\ref{eq:tmp3}) we deduce
$$
\sum_{k=1}^{n-1} (-1)^k p^k \sin(kz) = \frac{-p\sin(z)-(-1)^{n}p^{n}(p\sin((n-1)z)+\sin(nz)) }{1+2p\cos(z)+p^2} 
$$
Thus for any $M \ge 0$, $z \in (0,\pi)$ and $x \in (0,\beta)$
\begin{align*}
& \left| \sum_{m=0}^M   (-1)^{m}\sin((m+1)z) \left(\frac{x}{\beta}\right)^{m}\right| \\
& =  \left|\beta^{2}  
\frac{\big(\frac{x}{\beta}\big)^{M+1}(-1)^M \left(\frac{x}{\beta}\sin(z(M+1)) +\sin(z(M+2))\right) +\sin(z)}{x^2 + 2x\beta\cos(z) + \beta^2} \right|
< \frac{3}{\sin(z)^2}\,.
\end{align*}
Hence by  dominated convergence theorem we get
\begin{align*}
&\alpha \sin{\rho\pi}\int_0^\beta \frac{x^\alpha}{1+x^\alpha}\frac{1}{\beta^2 + 2x\beta\cos(\rho\pi)+x^2}dx \\
&= \alpha \sin{\rho\pi} \int_0^\beta \frac{x^\alpha}{1+x^\alpha}\frac{1}{\beta^2((\frac{x}{\beta})^2 + 2\frac{x}{\beta}\cos(\rho\pi)+1)}dx\\
&= \alpha \int_0^\beta \frac{x^\alpha}{1+x^\alpha} \sum_{m=0}^\infty \beta^{-2} (-1)^{m}\sin(\rho(m+1)\pi) \left(\frac{x}{\beta}\right)^{m}dx \\
&= \alpha \int_0^\beta  \lim_{M \to \infty} \sum_{m=0}^M \frac{x^\alpha}{1+x^\alpha} \beta^{-2} (-1)^{m}\sin(\rho(m+1)\pi) \left(\frac{x}{\beta}\right)^{m}dx \\
&=\sum_{m=0}^\infty  (-1)^{m} \beta^{-2-m}\sin(\rho(m+1)\pi)  \int_0^\beta \frac{ \alpha  x^{\alpha+m}}{1+x^\alpha} dx\,.
\end{align*}
By Lemma \ref{lem:int0b} 
\begin{align*}
\int_0^\beta \frac{ \alpha  x^{\alpha+m}}{1+x^\alpha} dx  = \int_0^{\beta^\alpha} \frac{ y^{(m+1)/\alpha}}{1+y} dy 
 = \sum_{k =0}^\infty \frac{(-1)^{k}\beta^{\alpha(k+1)+m+1}}{k+1+(m+1)/\alpha}\,.
\end{align*}
Consequently 
\begin{align}
&\alpha \sin{\rho\pi} \int_0^\beta \frac{x^\alpha}{1+x^\alpha}\frac{1}{\beta^2 +2x\beta\cos(\rho\pi)+x^2}dx \nonumber\\
&= \alpha\sum_{m=0}^\infty  \sum_{k =0}^\infty \frac{(-1)^{k+m}\beta^{\alpha(k+1)-1}\sin(\rho(m+1)\pi) }{\alpha(k+1)+(m+1)}\,.\label{eq:int0b}
\end{align}
Now we calculate $\int_\beta^\infty$. Similarly by the dominated convergence theorem we get
\begin{align*}
&\alpha \sin(\rho\pi) \int_\beta^\infty \frac{x^\alpha}{1+x^\alpha}\frac{1}{\beta^2 +2x\beta\cos(\rho\pi)+x^2}dx\\
 &= \alpha \sin(\rho\pi) \int_\beta^\infty \frac{x^\alpha}{1+x^\alpha}\frac{1}{x^2(1 +2\frac{\beta}{x}\cos(\rho\pi)+ (\frac{\beta}{x})^2)}dx\\
&= \alpha \int_\beta^\infty \frac{x^{\alpha-2}}{1+x^\alpha} \sum_{m=0}^\infty (-1)^{m} \left(\frac{\beta}{x}\right)^{m}\sin(\rho(m+1)\pi) dx \\
&=\sum_{m=0}^\infty  (-1)^{m} \beta^{m} \sin(\rho(m+1)\pi) \int_\beta^\infty \frac{ \alpha  x^{\alpha-2-m}}{1+x^\alpha} dx\\
&=\sum_{m=0}^\infty  (-1)^{m} \beta^{m} \sin(\rho(m+1)\pi) \int_{\beta^\alpha}^\infty \frac{ y^{-(1+m)/\alpha}}{1+y} dy\,.
\end{align*}
Since $\alpha \in \mathcal{A}$  by Lemma \ref{lem:intbinfty} we get
$$
\int_{\beta^\alpha}^\infty \frac{ y^{-(1+m)/\alpha}}{1+y} dy = \frac{\pi}{\sin(\frac{m+1}{\alpha}\pi)} - 
\sum_{k=0}^\infty \frac{(-1)^{k}\beta^{\alpha(k+1)-(m+1)}}{k+1 -(m+1)/\alpha}\,.
$$
Therefore by Lemma \ref{lem:Liouv1}
\begin{align}
&\alpha \sin{\rho\pi}\int_\beta^\infty \frac{x^\alpha}{1+x^\alpha}\frac{1}{\beta^2 +2x\beta\cos(\rho\pi)+x^2}dx \nonumber\\
&=\pi \sum_{m=0}^\infty  \frac{ (-1)^{m} \beta^{m} \sin(\rho(m+1)\pi)}{\sin(\frac{m+1}{\alpha}\pi)} \label{eq:intbinfty} \\
&- \alpha\sum_{m=0}^\infty  \sum_{k =0}^\infty \frac{(-1)^{k+m}\beta^{\alpha(k+1)-1}\sin(\rho(m+1)\pi)}{\alpha(k+1)-(m+1)}\,. \nonumber
\end{align}
Hence by (\ref{eq:int0b}), (\ref{eq:intbinfty}) and (\ref{eq:tmp1}) we get
\begin{align*}
&\frac{1}{\pi}\int_0^\infty \frac{x^\alpha}{1+x^\alpha}\frac{\alpha \sin{\rho\pi}}{\beta^2 +2x\beta\cos(\rho\pi)+x^2}dx 
= \sum_{m=0}^\infty  \frac{ (-1)^{m} \beta^{m} \sin(\rho(m+1)\pi)}{\sin(\frac{m+1}{\alpha}\pi)} \\
& + \frac{2\alpha}{\pi} \sum_{m=0}^\infty  \sum_{k =0}^\infty \frac{(-1)^{k+m}\beta^{\alpha(k+1)-1}(m+1)\sin(\rho(m+1)\pi)}{(m+1)^2 - (\alpha(k+1))^2} \\
&= \sum_{m=0}^\infty  \frac{ (-1)^{m} \beta^{m} \sin(\rho(m+1)\pi)}{\sin(\frac{m+1}{\alpha}\pi)} + \alpha \sum_{k =0}^\infty \frac{(-1)^{k}\beta^{\alpha(k+1)-1}\sin(\alpha\rho(k+1)\pi)}{\sin(\alpha(k+1)\pi)}. 
\end{align*}
The change of order of summation in the second line is justified by Lemma \ref{lem:Liouv2} and the Lebesgue theorem. 
\end{proof}
Now Theorem \ref{thm:main} follows easily from Lemma \ref{lem:main}. 

\section{Remarks and applications}

\begin{remark}\label{rem:1}\rm
Put 
\begin{equation}
g_k(a,x) = \sum_{m=1}^\infty\frac{ x^{m}U_{k-1}(\cos(m\pi a))}{m}\,,
\end{equation}
where $U_k(x)$ are the Chebyshev polynomials of the second type (we put $U_{-1} \equiv 0$). 
If $\rho +k=  l/\alpha$ (like in \cite{MR905336}), $l \ge 0$ and $k \ge 1$ we obtain for all $\alpha \in (0,2]$
\begin{equation}\label{eq:Doneycase}
g(\beta) = g_k(\alpha, (-1)^{l+1} \beta^\alpha ) - g_l(1/\alpha, (-1)^{k+1} \beta)\,,
\end{equation}
We note that sums above correspond to the function $f_k$ defined in \cite{MR905336}.
\end{remark}
\begin{proof}
First suppose $\alpha \in \mathcal{A}$ and $\rho = l/\alpha-k$. Since $U_k(\cos(x)) = \frac{\sin((k+1)x)}{\sin x}$ we get (\ref{eq:Doneycase}) for all $\alpha \in \mathcal{A}$. Now for $\alpha \in (0,2] \setminus \mathcal{A}$ we take $\mathcal{A} \ni \alpha_n \to \alpha$ and $\rho_n =  l/\alpha_n-k$. Passing to the limit we get (\ref{eq:Doneycase}) for $\alpha \in (0,2]$.  
\end{proof}
Using formulas (\cite[1.342.4, 1.342.2, 1.448.2]{MR2360010}) 
$$
U_{k-1}(\cos(z)) = 
\begin{cases} 
2 \sum\limits_{n=0}^m \cos((2n+1)z) & \mbox{for $k=2m+2$}\,,\\
1+ 2 \sum\limits_{n=1}^m \cos(2nz) & \mbox{for $k=2m+1$}\,,\\
\end{cases}
$$
$$
2 \sum_{m=1}^\infty\frac{x^{m}\cos(m z)}{m} = - \log(x^2 - 2x\cos(z) +1)\,,
$$
the functions $g_k(a,x)$ may be expressed by finite sums
$$
-g_k(a,x) =
\begin{cases}
\sum\limits_{n=0}^{k/2-1} \log(x^2 - 2x\cos((2n+1) a \pi) +1) & \mbox{for even $k$} \,, \\
\log(1-x) +\sum\limits_{n=1}^{(k-1)/2} \log(x^2 - 2x\cos(2n a \pi) +1) & \mbox{for odd $k$} \,.
\end{cases}
$$

\begin{example}\rm
Let $k=l=1$ then $\alpha \in (0,1)$ and 
$$
g(\beta) = - \sum_{m=1}^\infty \beta^m / m + \sum_{m=1}^\infty \beta^{\alpha m} / m = -\log(1 - \beta^\alpha) + \log(1-\beta)\,.
$$
Hence $\kappa(1,\beta) = \tilde{C}\frac{1-\beta}{1-\beta^\alpha }$.
\end{example}
A first application of Theorem 1 is to obtain new expressions for the functions $g'(\beta),g(\beta)$
and consequently $\kappa(1,\beta)$ and $ \kappa(\gamma,\beta)$ for the values of $\beta$
not concerned by the results of \cite{MR905336}.
\begin{prop}\label{lem:rational}
Let $\alpha \in \QQ \cap (0,2]$ and $\beta \in (0,1)$. Then 
\begin{align}
 g'(\beta)  &=  \sum_{\genfrac{}{}{0pt}{}{m=1,}{\frac{m}{\alpha} \not\in \NN}}^\infty \frac{(-1)^{m+1} \beta^{m-1}\sin(\rho m\pi)}{\sin(\frac{m\pi}{\alpha})} + \alpha\sum_{\genfrac{}{}{0pt}{}{k=1,}{\alpha k \not\in \NN}}^\infty \frac{(-1)^{k+1} \beta^{\alpha k-1}\sin(\rho\alpha k \pi)}{\sin(\alpha k\pi)} \nonumber\\ 
&  + \frac{\alpha\log(\beta)}{\pi}\sum_{\genfrac{}{}{0pt}{}{m=1,}{\frac{m}{\alpha} \in \NN}}^\infty (-1)^{m+\frac{m}{\alpha}} \beta^{m-1}\sin(\rho m\pi)  \label{eq:rational}\\
&+  \alpha\rho \sum_{\genfrac{}{}{0pt}{}{k=1,}{\alpha k \in \NN}}^\infty (-1)^{k(\alpha+1)}\beta^{\alpha k-1} \cos(\alpha\rho k\pi)\,. \nonumber
\end{align}
\end{prop}
\begin{proof}
Let $\alpha = \frac{p}{q}$. Like in Remark \ref{rem:1} we take $\mathcal{A} \ni \alpha_j = \frac{p}{q} + \frac{\sqrt{2}}{j}$. We obtain result by passing to the limit $j \to \infty$ in the expression
\begin{align*}
\sum_{\genfrac{}{}{0pt}{}{m=1,}{\frac{m}{\alpha} \not\in \NN}}^\infty \frac{(-1)^{m+1} \beta^{m-1}\sin(\rho m\pi)}{\sin(\frac{m\pi}{\alpha_j})} + \alpha_j\sum_{\genfrac{}{}{0pt}{}{k=1,}{\alpha k \not\in \NN}}^\infty \frac{(-1)^{k+1} \beta^{\alpha_j k-1}\sin(\rho\alpha_j k \pi)}{\sin(\alpha_j k\pi)}\\
+ \sum_{\genfrac{}{}{0pt}{}{m=1,}{\frac{m}{\alpha}  \in \NN}}^\infty \frac{(-1)^{m+1} \beta^{m-1}\sin(\rho m\pi)}{\sin(\frac{m\pi}{\alpha_j})} + 
\alpha_j\sum_{\genfrac{}{}{0pt}{}{k=1,}{\alpha k \in \NN}}^\infty \frac{(-1)^{k+1} \beta^{\alpha_j k-1}\sin(\rho\alpha_j k \pi)}{\sin(\alpha_j k\pi)}\\
\end{align*}
By Lemma \ref{lem:Liouv1} we pass with limit under sum signs. The first two terms obviously converge to the first two terms in (\ref{eq:rational}).  If we take $m=np$, $k=nq$, the second line is equal to
\begin{align}
&\sum_{n=1}^\infty \left(\frac{(-1)^{np+1} \beta^{n p-1} \sin(\rho n p \pi)}{\sin(n p \pi/\alpha_j)} + \alpha _j\frac{(-1)^{nq+1}  \beta^{n q \alpha_j-1} \sin(\rho n q \alpha_j \pi)}{\sin(n q \alpha_j \pi)} \right) \nonumber\\
& \stackrel{j \to \infty}{\longrightarrow}  
\sum_{n=1}^\infty \frac{(-1)^{n q+1}(-\beta)^{n p -1}  p (\pi \rho \cos(n p \pi \rho) + \log(\beta) \sin(n p \pi \rho))}{\pi q} \label{eq:rat2}
\end{align}
and the assertion of the  proposition holds. For the detailed proof of (\ref{eq:rat2}) we refer to Appendix.
\end{proof}
\begin{remark}\rm
In fact Proposition \ref{lem:rational} holds for $\alpha \in \mathcal{A} \cup (\QQ \cap 
(0,2])$  and (\ref{eq:rational}) may be treated as a generalization of Lemma \ref{lem:main}. 
\end{remark}
\begin{example} \rm
Let $\alpha = 1/2$. Then $p=1$ and $q=2$. We get
\begin{align*}
 g'(\beta) &= \frac{1}{2}\sum_{k=0}^\infty (-1)^k \beta^{k-\frac{1}{2}}\sin(\rho (k + \tfrac{1}{2}) \pi)\\
& + \frac{1}{2\pi}\sum_{n=1}^\infty (-1)^n \beta^{n -1} (\rho\pi \cos(n \rho \pi) + \log(\beta) \sin(n \rho \pi ))\\
& = \frac{\frac{(1 + \beta) \cos((\pi \rho)/2)}{2\sqrt{\beta}} - \frac{\rho(\beta +\cos(\pi \rho))}{2} - \frac{\log(\beta)\sin(\pi \rho)}{\pi}}{\beta^2+2\beta\cos(\pi \rho) +1}
\end{align*}
Analogous simple expressions can be given for other rational $\alpha$ not  covered by the results
of  \cite{MR905336}.
\end{example}
Further applications of formula (\ref{eq:thm}) from Theorem 1  are planned in the  forthcoming paper
\cite{PGTJ} where symmetric $\alpha$-stable processes $X_t$  in $\RR$ are considered.  The starting point is the formula
(see \cite{MR2250061})
\begin{equation}\label{eq:FactorId}
\int_0^\infty\int_0^\infty e^{-\eta t}e^{-\theta x}\EE_x(e^{-\gamma X_t};\; \tau>t)\,dt\,dx = \frac{1}{(\theta+\gamma)\,\kappa(\eta,\gamma)\,\kappa(\eta,\theta)}\,,
\end{equation}
where $\tau = \tau_{(0,\infty)}$ is the first exit time from $(0,\infty)$ of the process $X_t$.

A better knowledge of $\kappa$ then permits to get from (\ref{eq:FactorId})
more  information about the law of $\tau$.

\section{Appendix}
Here we give a detailed proof of (\ref{eq:rat2}).
\begin{lem} Let $\alpha_j=\frac{p}{q} +\frac{\sqrt{2}}{j}$. We have

\begin{align}
&\sum_{n=1}^\infty \left(\frac{(-1)^{np+1} \beta^{n p-1} \sin(\rho n p \pi)}{\sin(n p \pi/\alpha_j)} + \alpha _j\frac{(-1)^{nq+1}  \beta^{n q \alpha_j-1} \sin(\rho n q \alpha_j \pi)}{\sin(n q \alpha_j \pi)} \right) \nonumber\\
& \stackrel{j \to \infty}{\longrightarrow}  
\sum_{n=1}^\infty \frac{(-1)^{n q+1}(-\beta)^{n p -1}  p (\pi \rho \cos(n p \pi \rho) + \log(\beta) \sin(n p \pi \rho))}{\pi q} 
\end{align}
\end{lem}
\begin{proof}
Let us call 
$$
F_1(n,j)= \frac{(-1)^{np+1} \sin(\rho n p \pi)}{\sin(n p \pi/\alpha_j)} +
\alpha _j\frac{(-1)^{nq+1}   \sin(\rho n p \pi)}{\sin(n q \alpha_j \pi)} 
$$
$$
F_2(n,j)= \alpha _j\frac{(-1)^{nq+1} ( \sin(\rho n q \alpha_j \pi)-\sin(\rho n p \pi))}{\sin(n q \alpha_j \pi)}
$$
$$
F_3(n,j)= \alpha _j\frac{(-1)^{nq+1}   \sin(\rho n q \alpha_j \pi)}{\sin(n q \alpha_j \pi)}(\beta^{nq\alpha_j-np}-1) 
$$
The proof of Lemma \ref{eq:rat2} consists of two parts:

1) {\it Term by term convergence:} 
We note that 
$$\sin(\tfrac{n p \pi}{\alpha_j}) = (-1)^{nq+1} \sin(\tfrac{nq^2\sqrt{2} \pi}{pj+\sqrt{2}})\,, \qquad \sin(n q \alpha_j \pi) = (-1)^{np} \sin(\tfrac{nq\sqrt{2}\pi}{j})\,.$$ 
Hence for fixed $n$ and large $j$ we have
\begin{align*}
\left|\frac{F_1(n,j)}{\sin(\rho n p \pi)}\right| &= \left|\frac{\sin(\frac{nq\sqrt{2}\pi}{j}) - \frac{pj + \sqrt{2}q}{qj} \sin(\frac{nq^2\sqrt{2} \pi}{pj+\sqrt{2}q})}{\sin(\frac{nq\sqrt{2}\pi}{j})\sin(\frac{nq^2\sqrt{2} \pi}{pj+\sqrt{2}q})}\right|\\
& \le \frac{\sum\limits_{k=1}^\infty \frac{1}{(2k+1)!}\left(\left(\frac{nq\sqrt{2}\pi}{j}\right)^{2k+1} + \frac{2p}{q} \left(\frac{nq^2\sqrt{2} \pi}{pj+\sqrt{2}q}\right)^{2k+1} \right)}{\frac{nq\sqrt{2}\pi}{2j}\frac{nq^2\sqrt{2} \pi}{2(pj+\sqrt{2}q)}}
 \le \frac{Kn}{j} \stackrel{j \to \infty}{\longrightarrow} 0\,,
\end{align*}
where $K$ is some constant independent of $n$ and $j$. Therefore $\lim\limits_{j \to \infty} F_1(n,j)=0$. Further
\begin{align*}
\lim\limits_{j \to \infty} F_2(n,j) &= \lim\limits_{j \to \infty} 2 \alpha_j (-1)^{n(p+q)+1} \frac{\sin(\frac{\rho n q \sqrt{2}\pi}{2 j}) \cos(\rho n \pi(q \alpha_j+p)/2)}{\sin(nq\sqrt{2}\pi/j)}\\
& = \frac{(-1)^{n ( p+q)+1}  p  \rho \cos( \rho n p \pi )} {q}\,.
\end{align*}
Similarly
\begin{align*}
\lim\limits_{j \to \infty} F_3(n,j) &= \lim\limits_{j \to \infty} \alpha _j (-1)^{n(p+q)+1}   \sin(\rho n q \alpha_j \pi) \frac{(\beta^{nq\sqrt{2}/j}-1)}{\sin(n q \sqrt{2} \pi/j)}\\
&= \frac{(-1)^{n (p+q)+1} p \log(\beta) \sin( n p \pi\rho )}{\pi q}
\end{align*}

2) {\it Uniform integrability with respect to the measure $\mu=\sum_{n=1}^\infty \beta^{np-1}\delta_n$}: We will show that for each $k =1,2,3$ we have
$$
\sup_{j \in \NN} \sum_{n=1}^\infty  |F_k(n,j)| \beta^{np-1} < \infty\,,
$$
and for every $\varepsilon >0$ there exists $\delta>0$ such that
$$
\sup_{j \in \NN} \sum_{n \in G}  |F_k(n,j)| \beta^{np-1} < \varepsilon\,, 
$$
whenever $\mu(G) <\delta$.

From part 1) of the proof we see that for $n < j/(2q\sqrt{2})$, we have $F_k(n,j) < C$, where $C$ does not depend on $n,j$ and $k=1,2,3$. Denote $G_j = \{m \in \NN \colon m < j/(2q\sqrt{2})$. Let $k \in \NN$ be the closest integer to $nq\sqrt{2}/j$. Then by diophantine approximation
\begin{align*}
|\sin(nq\alpha_j\pi)| &= |\sin((k-nq\sqrt{2}/j)\pi)| \ge \frac{|k-nq\sqrt{2}/j|}{2} \\
& = \frac{nq}{2j}\left|\sqrt{2} -\frac{kj}{nq}\right| \ge \frac{nq}{j}\frac{c_1}{(nq)^2} = \frac{c_2}{nj}\,.
\end{align*}
Similarly we show that
$$
|\sin(n p \pi/\alpha_j)| = \left|\sin\left(\frac{nq^2\sqrt{2} \pi}{pj+\sqrt{2}q}\right)\right| \ge \frac{c_3}{nj}\,.
$$
Therefore
\begin{align*}
\sup_{j \in \NN}\sum_{n=1}^\infty |F_1(n,j)| \beta^{np-1} \le \sup_{j \in \NN} \left(\sum_{n \in G_j} C \beta^{np-1} +  c\sum_{n \in \NN \setminus G_j} nj \beta^{np-1}\right) < \infty\,. 
\end{align*}
Now let $\varepsilon >0$. First we note that 
$$
b_j = \sum_{n \in \NN \setminus G_j} n j \beta^{np-1} \to 0\,, \qquad {\rm if \,} j \to \infty\,.  
$$
Hence there is $j_0 \in \NN$ such that for $j > j_0$ we have $b_j < \varepsilon/3$.
We take $n_0 \in \NN$ such that $\sum_{n=n_0}^\infty n \beta^{np-1} < \varepsilon/(3cj_0)$ and $\sum_{n=n_0}^\infty \beta^{np-1} < \varepsilon/(3C)$. Now let $\delta = \beta^{n_0 p -1}$. If $\mu(G) < \delta$ then $G \subset \{n_0,n_0+1, \dots\}$ and 
\begin{align*}
\sup_{j \in \NN}\sum_{n \in G} |F_1(n,j)| \beta^{np-1} &\le \sup_{j \in \NN} \sum_{n \in G \cap G_j} C \beta^{np-1} +  \sup_{j \in \NN}  c\sum_{n \in G \setminus G_j} nj \beta^{np-1}\\ 
&\le C \sum_{n\in G} \beta^{np-1} + c \sup_{j > j_0} \sum_{n \in \NN \setminus G_j} n j \beta^{np-1} + c \sum_{n \in G} n j_0 \beta^{np-1}\\
& \le \frac{\varepsilon}{3} + \frac{\varepsilon}{3} + \frac{\varepsilon}{3} =\varepsilon\,.
\end{align*}
In the same way we prove uniform integrability of $F_2(n,j)$ and $F_3(n,j)$ and we obtain the assertion of the lemma.
\end{proof}
\section*{Acknowledgements}
We thank Zbigniew Palmowski for introducing us into the subject. We also thank Michel Waldschmidt for discussions about this paper.

\bibliographystyle{abbrv}

\end{document}